\documentclass[12pt]{amsart}
\usepackage{amssymb,amsmath,mathrsfs}
\usepackage{amsthm,amsfonts}
\usepackage{paralist}
\usepackage{color}
\usepackage{epsfig}
\usepackage{graphicx}
\usepackage{appendix}

\begin{document}

\def\Xint#1{\mathchoice
  {\XXint\displaystyle\textstyle{#1}}%
  {\XXint\textstyle\scriptstyle{#1}}%
  {\XXint\scriptstyle\scriptscriptstyle{#1}}%
  {\XXint\scriptscriptstyle\scriptscriptstyle{#1}}%
  \!\int}
\def\XXint#1#2#3{{\setbox0=\hbox{$#1{#2#3}{\int}$}
  \vcenter{\hbox{$#2#3$}}\kern-.5\wd0}}
\def\ddashint{\Xint=}
\def\dashint{\Xint-}

\newcommand{\R}{\mathbb{R}}

\def\pd#1#2{\frac{\partial#1}{\partial#2}}
\def\dfrac{\displaystyle\frac}
\let\oldsection\section
\renewcommand\section{\setcounter{equation}{0}\oldsection}
\renewcommand\thesection{\arabic{section}}
\renewcommand\theequation{\thesection.\arabic{equation}}

\newtheorem{thm}{Theorem}[section]
\newtheorem{lem}{Lemma}[section]
\newtheorem{definition}{Definition}[section]
\newtheorem{prop}{Proposition}[section]
\newtheorem{remark}{Remark}[section]
\newtheorem{cor}{Corollary}[section]
\allowdisplaybreaks

\title[]{Solutions of  Fully Nonlinear Nonlocal
Systems}

\author[Pengyan Wang and Mei Yu]{}

\email{wangpy119@126.com}
\email{yumei@nwpu.edu.cn}

\maketitle

\centerline{\scshape Pengyan Wang, Mei Yu}
\medskip
{\footnotesize\centerline{Department of Applied Mathematics , Northwestern Polytechnical University, } \centerline{Xi'an,  710129,
Shaan xi, P.R. China} }

\begin{abstract}
In this paper we consider   the system  involving fully nonlinear  nonlocal operators:
$$ \left\{
\begin{array}{ll}
F_{\alpha}(u(x)) = C_{n,\alpha} PV \int_{{R}^n} \frac{G(u(x)-u(y))}{|x-y|^{n+\alpha}} dy=f(v(x)),\\
F_{\beta}(v(x)) = C_{n,\beta} PV \int_{{R}^n} \frac{G(v(x)-v(y))}{|x-y|^{n+\beta}} dy=g(u(x)).
\end{array}
\right.
$$
  A \textit{narrow region principle} and   a \textit{decay at infinity} for the system  for carrying on  the method of moving planes are established. Then we
prove the radial symmetry and monotonicity for   positive
solutions to the nonlinear system  in the whole space.  Non-existence of positive
solutions to the nonlinear system on a half space is proved.
\end{abstract}

\section{Introduction}
In this paper, we consider the nonlinear system  involving fully nonlinear nonlocal operators: $$\left\{
\begin{array}{ll}
F_{\alpha}(u(x)) = f(v(x)),\\
F_{\beta}(v(x)) =g(u(x)),
\end{array}
\right.$$
with $$F_{\alpha}(u(x))
= C_{n,\alpha} PV \int_{{R}^n} \frac{G(u(x)-u(y))}{|x-y|^{n+\alpha}} dy,$$
where $PV$ stands for the Cauchy principal value, $G$ is at least local Lipschitz continuous,  $G(0)=0$ and $0<\alpha, \beta<2$. The operators $F_\alpha$ was introduced by
Caffarelli and Silvestre in \cite{CS1}.

In order the integral to make sense, we require  $$u \in C^{1,1}_{loc}\cap L_\alpha,~v\in C^{1,1}_{loc}\cap L_\beta$$
with $$ L_{\alpha }=\{u:R^{n}\rightarrow R \mid \int_{R^{n}} \frac{|u(x)|}{1+|x|^{n+\alpha}}dx <\infty   \} , $$
and $L_\beta$ having a similar meaning.

In the special case when $G(\cdot)$ is an identity map, $F_\alpha$ becomes the usual fractional Laplacian $(-\Delta)^{\frac{\alpha}{2}}.$  The   nonlocal nature of fractional operators makes them difficult to study. To circumvent this,   Caffarelli and Silvestre \cite{CS} introduced the \textit{extension method} which turns the nonlocal problem involving the fractional Laplacian into a local one in higher dimensions. This method has been applied successfully to treat equations involving the fractional Laplacian
and a series of fruitful results has been obtained (see \cite{BCPS}, \cite{CZ}, etc.). One can
also use \textit{the integral equations method}, such as \textit{the method of moving planes in integral forms} (see \cite{CC}, \cite{CD},  \cite{ZCCY}, \cite{LZ}, \cite{LZr}) and \textit{regularity lifting }to investigate equations involving fractional Laplacian by  showing that they are equivalent to corresponding integral equations (see \cite{CFY},  \cite{CLO}, \cite{CLO1} and the references therein).  For more articles concerning the method of moving planes for nonlocal equations and
for integral equations,   see \cite{FL}, \cite{HLZ}, \cite{HWY}, \cite{LLM}, \cite{LZ2}, \cite{MC}, \cite{MZ}, \cite{LZ} and the references therein.

For the fully nonlinear nonlocal equations, so far as we know, there is neither   any corresponding {\em extension method} nor equivalent integral equations that one can work at. A probable  reason is that very few results were obtained for  fully nonlinear nonlocal operator. In \cite{CLL}, Chen, Li and Li developed a new method that can deal directly with these   nonlocal operators. Inspired by the idea,  we extend the method in \cite{CLL} to fully nonlinear nonlocal systems and  consider the nonlinear systems  involving fully nonlinear nonlocal operators
 \begin{equation}\label{eq:a1}
 \left\{
\begin{array}{ll}
F_{\alpha}(u(x)) = f(v(x)),\\
F_{\beta}(v(x)) =g(u(x)), &\mbox~ x\in R^n,\\
u(x)>0,v(x)>0,&\mbox~ x\in R^n,
\end{array}
\right.\end{equation}
and
\begin{equation}\label{eq:a2}
\left\{
\begin{array}{ll}
F_{\alpha}(u(x)) = f(v(x)),\\
F_{\beta}(v(x)) =g(u(x)), &\mbox~ x\in R^n_+,\\
u(x)\equiv 0,v(x)\equiv 0,&\mbox~ x\not\in R^n_+,
\end{array}
\right.
\end{equation}
where $f$ and $g$  are nonnegative continuous and nondecreasing functions.

We first establish the  \textit{narrow region principle}  and \textit{decay at infinity} for  the system  which play important roles in carrying out the method of moving planes.
To state them, denote by
$$T_\lambda=\{x\in R^n| x_1=\lambda \}$$
 the moving plane, $$\Sigma_\lambda=\{x\in R^n|x_1<\lambda\}$$  the left region   of the plane $T_\lambda$,   $$x^\lambda=(2\lambda-x_1,x_2,\cdots, x_n)$$  the reflection
of $x$ about  $T_\lambda$, and denote $$U_\lambda(x)=u_\lambda(x)-u(x)~ \mbox{and}~
V_\lambda(x)=v_\lambda(x)-v(x).$$ For simplicity of notations,  we
stand for $U_\lambda(x)$ by $U(x)$ and $V_\lambda(x)$ by $V(x)$ in the sequel.

\begin{thm}\label{thma2}(Narrow Region Principle )
Let $\Omega$ be a bounded narrow region in $\Sigma_\lambda$  contained in
 $$\{x |  \lambda-l<x_1<\lambda  \}$$
 with small $l>0$.  Suppose  that $U(x)\in L_\alpha\cap C_{loc}^{1,1}(\Omega), V(x)\in L_\beta \cap C_{loc}^{1,1}(\Omega), $ and $U(x),V(x)$ are lower semi-continuous on $\bar{\Omega}$.
 If
 $c_i(x)\leq 0,i=1,2,$ are bounded from below in $\Omega$, $U(x)$ and $V(x)$ satisfy
\begin{equation}\label{eq:a3}
\left\{\begin{array}{ll}
F_\alpha (u_\lambda(x))-F_\alpha (u(x))+c_1(x)V(x)\geq0, \\
F_\beta (v_\lambda(x))-F_\beta (v(x))+c_2(x)U(x)\geq0, &\quad x \in \Omega,\\
U(x), V(x)\geq 0, &\quad x \in \Sigma_\lambda\backslash\Omega,\\
U(x^\lambda)=-U(x),
V(x^\lambda)=-V(x), &\quad x \in \Sigma_\lambda,
\end{array}
\right.
\end{equation}
then we have for sufficiently small $l$,
\begin{equation}\label{eq:a4}
U(x),V(x) \geq0 \mbox{ in } \Omega;\end{equation}
if $\Omega$ is unbounded, the conclusion still holds under the conditions
$$\underset{|x|\rightarrow \infty}{\underline{\lim}}U(x),\underset{|x|\rightarrow \infty}{\underline{\lim}}V(x)\geq0;$$
furthermore, if $U(x)$ or $V(x)$ attains 0 somewhere in $\Sigma_\lambda$, then
\begin{equation}\label{eq:a5}U(x)=V(x)\equiv 0,~ x\in R^n.\end{equation}

\end{thm}
We call \eqref{eq:a5} the strong maximum principle later. As we can see from the proof,  to ensure \eqref{eq:a5}, $\Omega$ does not need to be narrow.
\begin{thm}\label{thma1}( Decay at Infinity)
Let $\Omega$ be a bounded or unbounded  domain in $R^n$. Assume that $U(x)\in  C^{1,1}_{loc}(\Omega)\cap L_\alpha(R^n), V(x)\in C^{1,1}_{loc}(\Omega)\cap L_\beta(R^n),$ $U(x)$ and $V(x) $ are lower semi-continuous on $ \bar {\Omega}$. If $U(x)$ and $V(x) $ satisfy
\begin{equation}\label{eq:a11}
\left\{
\begin{array}{ll}
F_{\alpha}(u_\lambda(x))-F_{\alpha}(u(x)) +c_1(x)V(x)\geq 0,\\
F_{\beta}(v_\lambda(x))-F_{\beta}(v(x)) +c_2(x)U(x)\geq 0, &\mbox  x\in \Omega,\\
U(x),V(x)\geq0, &\mbox x \in {\Sigma_\lambda} \backslash{\Omega},\\
U(x^\lambda) =-U(x),\\
V(x^\lambda) =-V(x), &\mbox x \in \Sigma_\lambda ,
\end{array}
\right.\end{equation}
with \begin{equation}\label{eq:aa11}c_1(x)\sim o(\frac{1}{|x|^\alpha}),~c_2(x)\sim o(\frac{1}{|x|^\beta}),~~ \mbox{for} ~|x|~~\mbox{large}, \end{equation}
and $$c_i(x)\leq 0,~i=1,2,$$
then there exists a constant  $R_0>0$  depending only on $c_i(x)$  such that if
$$U(\tilde{x})=\underset{\Omega} {min}U(x) <0~ \mbox{and} ~ V(\bar{x})=\underset{\Omega} {min}V(x) <0,$$
  then \begin{equation}\label{eq:aaa11}|\tilde{x}|\leq R_0~\mbox{or} ~ |\bar x|\leq R_0.\end{equation}
\end{thm}
Based on Theorems \ref{thma2} and  \ref{thma1}, we apply
 the  \textit{method of moving planes} to obtain symmetry and monotonicity of positive solutions to  \eqref{eq:a1} in $R^n$,  as well as nonexistence of positive solutions to \eqref{eq:a2} on the half space.

\begin{thm}\label{thma3}
   Assume that $u(x)\in L_\alpha(R^n)\cap C_{loc}^{1,1}(R^n)~ \mbox{and}~v(x)\in L_\beta(R^n)\cap C_{loc}^{1,1}(R^n)$ are  positive solutions of system \eqref{eq:a1}. Suppose that for some $ \gamma,\tau>0,$
     \begin{equation}\label{eq:aa1} v(x)=o(\frac{1}{|x|^\gamma}),~u(x)=o(\frac{1}{|x|^\tau}),~ \mbox{as}~ |x|\rightarrow \infty,\end{equation} and \begin{equation}\label{eq:aa2}f'(s)\leq s^q, ~g'(t)\leq t^p,~ \mbox{with} ~~ q\gamma\geq\alpha,~ p\tau\geq\beta  .\end{equation}
      Then  $u(x)$ and $v(x)$ must be radially symmetric and monotone decreasing about some point $x_0$ in $R^n$.
 \end{thm}

\begin{thm}\label{thma4}Assume that $u(x)\in L_\alpha\cap C_{loc}^{1,1}(R^n_+),~v(x)\in L_\beta\cap C_{loc}^{1,1}(R^n_+)$ are  nonnegative solutions of   system \eqref{eq:a2}. Suppose \begin{equation}\label{eq:a28}\underset{|x|\rightarrow \infty}{\lim}u(x)=0,~\underset{|x|\rightarrow \infty}{\lim}v(x)=0,\end{equation}
$f(v),g(u)$ are Lipschitz continuous in the range of $ v(x),u(x)$ respectively,  and $f(0)=0,~g(0)=0.$
Then   $u(x)\equiv0,~v(x)\equiv0$.
\end{thm}

In section 2,  we   prove Theorems \ref{thma2} and  \ref{thma1} with a key ingredient \eqref{eq:aa20} below.  In section 3,  the  proofs of  Theorems \ref{thma3}  and \ref{thma4} are given by using the previous results and the method of moving planes.

\section{Proofs of Theorems \ref{thma2} and \ref{thma1} }
 ~~~~~Let $$F_{\alpha}(u(x))= C_{n,\alpha} PV \int_{\mathbb{R}^n} \frac{G(u(x)-u(y))}{|x-y|^{n+\alpha}} dy=C_{n,\alpha} \lim_{\epsilon \rightarrow 0} \int_{\mathbb{R}^n\setminus B_{\epsilon}(x)} \frac{G(u(x)-u(y))}{|x-y|^{n+\alpha}} dy
.$$
Throughout this and next section, we assume
\begin{equation}\label{eq:a33}G\in C^1(R),~G(0)=0, \mbox{and} ~G'(t)\geq c_0>0, ~\mbox{for} ~t \in R.\end{equation}

Using the simple maximum principle in \cite{CLLg}, we prove the following strong maximum principle.
\begin{lem}\label{lemaa1}

Let $\Omega$ be a bounded domain in $R^n$. Assume that $u(x)\in  C^{1,1}_{loc}(\Omega)\cap L_\alpha(R^n),$ is lower semi-continuous on $\bar \Omega$,
and satisfies \begin{equation}\label{eq:a30}
 \left\{
\begin{array}{ll}
F_{\alpha}(u(x)) \geq 0, &\mbox x\in \Omega,\\
u(x)\geq0,&\mbox x\in \Omega^c.
\end{array}
\right.\end{equation}
If $u(x)$  attains 0 somewhere in $\Sigma_\lambda$, then
$$u(x)\equiv 0, ~x\in R^n.$$
\end{lem}
{\bf Proof .}
 If $u(x)$ is not identical to 0, there exists
an $x^0$ such that $u(x^0)=0$ and
$$\aligned F_{\alpha}(u(x))&=\int_{R^n} \frac{G(u(x^0)-u(z))}{|x^0-z|^{n+\alpha}}dz\\
&=\int_{R^n} \frac{G'(\Psi(z))[u(x^0)-u(z)]}{|x^0-z|^{n+\alpha}}dz\\
&\leq c_0 \int_{R^n} \frac{-u(z)}{|x^0-z|^{n+\alpha}}dz\\
&<0.
\endaligned$$
 This contradicts \eqref{eq:a30}  and the proof is ended.

{\bf Proof of Theorem \ref{thma2}.}

If \eqref{eq:a4} does not hold, without loss of generality, we assume $U(x)<0$ at some point in $\Omega$; then the lower semi-continuity of $U(x)$ on $\bar \Omega$ guarantees that there exists some $\tilde{x}\in \Omega $ such that
$$U(\tilde{x})=\underset{\Omega} {\min}U(x) <0.$$
And it deduces from the condition \eqref{eq:a3} that $\tilde{x}$ is in the interior of $\Omega.$
By the defining integral, we have
\begin{eqnarray}\label{eq:aa20}
 F_\alpha (u_\lambda(\tilde x))-F_\alpha (u(\tilde x))&=&C_{n,\alpha} PV\int_{R^{n}}\frac{G(u_\lambda (\tilde x)-u_\lambda (y))-G(u (\tilde x)-u(y))}{|\tilde x -y|^{n+\alpha}}dy \nonumber\\
 &=&C_{n,\alpha} PV\int_{\Sigma_\lambda}\frac{G(u_\lambda (\tilde x)-u_\lambda (y))-G(u (\tilde x)-u(y))}{|\tilde x -y|^{n+\alpha}}dy\nonumber\\&+&C_{n,\alpha} PV\int_{\Sigma_\lambda}\frac{G(u_\lambda (\tilde x)-u (y))-G(u (\tilde x)-u_\lambda(y))}{|\tilde x -y^\lambda|^{n+\alpha}}dy\nonumber\\
 &\leq& C_{n,\alpha} PV\int_{\Sigma_\lambda}\frac{G(u_\lambda (\tilde x)-u_\lambda (y))-G(u (\tilde x)-u(y))}{|\tilde x -y^\lambda|^{n+\alpha}}dy\nonumber\\&+&C_{n,\alpha} PV\int_{\Sigma_\lambda}\frac{G(u_\lambda (\tilde x)-u (y))-G(u (\tilde x)-u_\lambda(y))}{|\tilde x -y^\lambda|^{n+\alpha}}dy\nonumber\\
 &=&C_{n,\alpha} PV\int_{\Sigma_\lambda}\frac{2G'(\cdot)U(\tilde x)}{|\tilde x -y^\lambda|^{n+\alpha}}dy\nonumber\\
 &\leq& 2C_{n,\alpha} c_0 U(\tilde x)\int_{\Sigma_\lambda}\frac{1}{|\tilde x -y^\lambda|^{n+\alpha}}dy.
\end{eqnarray}
Let $D=B_{2l}(\tilde x)\cap \tilde \Sigma_\lambda $,   then
\begin{equation}\label{eq:a6}\aligned
\int_{\Sigma_\lambda}\frac{1}{|\tilde x -y^\lambda|^{n+\alpha}}dy
&\geq \int_{D}\frac{1}{|\tilde x -y|^{n+\alpha}}dy\\
&\geq\frac{1}{10}\int_{B_{2l}(\tilde x)}\frac{1}{|\tilde x -y|^{n+\alpha}}dy\\
&\geq\frac{1}{l^\alpha}.\endaligned
\end{equation}
Thus from \eqref{eq:aa20},
\begin{equation}\label{eq:a7}
F_\alpha (u_\lambda(\tilde x))-F_\alpha (u(\tilde x))\leq \frac{C U(\tilde x)}{l^\alpha}<0.
\end{equation}
Together \eqref{eq:a7} with \eqref{eq:a3}, it yields
\begin{equation}\label{eq:a8}
U(\tilde x)\geq -Cc_1(\tilde x)l^\alpha V(\tilde x)~\mbox{and}~V(\tilde x)\leq 0.
\end{equation}
We know from   \eqref{eq:a8} that there exists  $\bar x$ such that
$$V(\bar x)=\underset{\Omega}{\min} V(x)<0.$$
Similarly to \eqref{eq:a7}, it derives that
$$F_\beta (v_\lambda(\bar x))-F_\beta (v(\bar x))\leq \frac{C V(\bar  x)}{l^\beta}<0.$$
Combining it with \eqref{eq:a8}, we have for $l$ sufficiently small,
$$\aligned
0&\leq F_\beta (v_\lambda(\bar x))-F_\beta (v(\bar x))+c_2(\bar x)U(\bar x)\\
&\leq\frac{CV(\bar x)}{l^\beta}+c_2(\bar x)U(\tilde x)\\
&\leq C(\frac{V(\bar x)}{l^\beta}-c_2(\bar x)c_1(\tilde x)l^\alpha V(\tilde x))\\
&\leq C(\frac{V(\bar x)}{l^\beta}-c_2(\bar x)c_1(\tilde x)l^\alpha V(\bar x))\\
&\leq C \frac{V(\bar x)}{l^\beta}(1-c_1(\tilde x)c_2(\bar x)l^{\alpha+\beta})\\
&<0.
\endaligned$$
This contradiction shows that \eqref{eq:a4} must be true.

Next we prove \eqref{eq:a5}. Without loss of generality,  let us suppose that there exists $\eta\in \Omega$ such that
$$U(\eta)=0.$$
Then we use $\frac{1}{|x-y|}>\frac{1}{|x-y^\lambda|}$,~for $x,y\in \Sigma_\lambda$, to have
 \begin{equation}\label{eq:aa10}\aligned
 &F_\alpha (u_\lambda(\eta))-F_\alpha (u(\eta))\\&=C_{n,\alpha} PV\int_{R^{n}}\frac{G(u_\lambda (\eta)-u_\lambda (y))-G(u (\eta)-u(y))}{|\eta -y|^{n+\alpha}}dy\\
 &=C_{n,\alpha} PV\int_{\Sigma_\lambda}\frac{G(u_\lambda (\eta)-u_\lambda (y))-G(u (\eta)-u(y))}{|\eta -y|^{n+\alpha}}dy\\&+C_{n,\alpha} PV\int_{\Sigma_\lambda}\frac{G(u_\lambda (\eta)-u (y))-G(u (\eta)-u_\lambda(y))}{|\eta -y^\lambda|^{n+\alpha}}dy \\
 &=C_{n,\alpha} PV\int_{\Sigma_\lambda}[G(u_\lambda (\eta)-u_\lambda (y))-G(u (\eta)-u(y))](\frac{1}{|\eta -y|^{n+\alpha}}-\frac{1}{|\eta -y^\lambda|^{n+\alpha}})dy\\
 &+C_{n,\alpha} PV\int_{\Sigma_\lambda}\frac{G(u_\lambda (\eta)-u (y))-G(u (\eta)-u_\lambda(y))+G(u_\lambda (\eta)-u_\lambda (y))-G(u (\eta)-u(y))}{|\eta -y^\lambda|^{n+\alpha}}dy\\
 & =C_{n,\alpha} G'(\cdot)\int_{\Sigma_\lambda} (U(\eta)-U(y))(\frac{1}{|\eta -y|^{n+\alpha}}-\frac{1}{|\eta -y^\lambda|^{n+\alpha}})dy\\&+C_{n,\alpha} G'(\cdot)\int_{\Sigma_\lambda}\frac{2U(\eta)}{|\eta -y^\lambda|^{n+\alpha}}dy\\
&\leq -Cc_0 \int_{\Sigma_\lambda}U(y)(\frac{1}{|\eta -y|^{n+\alpha}}-\frac{1}{|\eta -y^\lambda|^{n+\alpha}})dy.
 \endaligned
\end{equation}
If $U(x)\not \equiv0$, then \eqref{eq:aa10} implies  $$F_\alpha (u_\lambda(\eta))-F_\alpha (u(\eta))<0.$$
Using it with \eqref{eq:a3}, it shows  $V(\eta)<0.$
This is a contradiction with \eqref{eq:a4}. Hence $U(x)$ must be identically $0$ in $\Sigma_\lambda$.
Since $$U(x^\lambda)=-U(x),~x\in \Sigma_\lambda,$$ it gives $$U(x)\equiv0,~x\in R^n.$$
Again from \eqref{eq:a3}, we see $$V(x)\leq 0,~x\in \Sigma_\lambda.$$
Since we already know  $$V(x)\geq 0,~x\in \Sigma_\lambda,$$
it must hold  $$V(x)=0,~x\in \Sigma_\lambda.$$
Recalling $V(x^\lambda)=-V(x),$  we arrive at
$$V(x)\equiv0,~x\in R^n.$$
Similarly, one can show that if $U(x)$ or $V(x)$ attains 0 at one point in $\Sigma_\lambda$, then both $U(x)$ and $V(x)$ are identically 0 in $R^n$. This completes the proof.

{\bf Proof of Theorem \ref{thma1}.} There exists  $\tilde x \in \Omega,$
such that $$U(\tilde x)=\underset{\Omega}{\min}U(x)<0.$$
Using \eqref{eq:aa20},  we have
$$
\aligned F_{\alpha}(u_\lambda(\tilde x)) -F_{\alpha}(u(\tilde x))
=&C_{n,\alpha}PV \int_{\Sigma_\lambda} \frac{G(u_\lambda(\tilde x)-u_\lambda(y))-G(u(\tilde x)-u(y))}{|\tilde x-y|^{n+\alpha}}\\+&C_{n,\alpha}PV \int_{ \Sigma_\lambda} \frac{G(u_\lambda(\tilde x)-u(y))-G(u(\tilde x)-u_\lambda(y))}{|\tilde x-y^\lambda|^{n+\alpha}}\\
\leq&C_{n,\alpha}PV \int_{\Sigma_\lambda}\frac{G'(\cdot)2U(\tilde x)}{|\tilde x-y^\lambda|^{n+\alpha}}dy\\
\leq&2C_{n,\alpha}c_0U(\tilde x)\int_{\Sigma_\lambda}\frac{1}{|\tilde x-y^\lambda|^{n+\alpha}}dy.
\endaligned
$$
For each fixed $\lambda$, there exists  $C>0$  such that for $\tilde x \in \Sigma_\lambda$ and $|\tilde x|$ sufficiently large,
\begin{equation}\label{eq:a9}
\int_{\Sigma_\lambda}\frac{1}{|\tilde x-y^\lambda|^{n+\alpha}}dy\geq \int_{B_{3|\tilde x|}(\tilde x)\backslash B_{2|\tilde x|}(\tilde x)} \frac{1}{|\tilde x-y|^{n+\alpha}}dy \sim \frac{C}{|\tilde x|^\alpha}.
\end{equation}
Hence
\begin{equation}\label{eq:a10}
F_{\alpha}(u_\lambda(\tilde x)) -F_{\alpha}(u(\tilde x)) \leq \frac{CU(\tilde x)}{|\tilde x|^\alpha}<0.
\end{equation}
Together \eqref{eq:a10} with \eqref{eq:a11}, it is easy to deduce
 \begin{equation}\label{eq:a12}
 V(\tilde x)<0,
 \end{equation}
and
\begin{equation}\label{eq:a13}
 U(\tilde x)\geq -Cc_1(\tilde x)|\tilde x|^\alpha V(\tilde x).
 \end{equation}
From \eqref{eq:a12}, there exists  $\bar x$ such that
$$V(\bar x)=\underset{\Omega}{\min}V(x)<0.$$
Similarly to \eqref{eq:a9}, we can derive
\begin{equation}\label{eq:aa13}F_{\beta}(v_\lambda(\bar x)) -F_{\beta}(v(\bar x)) \leq \frac{CV(\bar x)}{|\bar x|^\beta}<0.\end{equation}
Combing \eqref{eq:a11} and \eqref{eq:a13}, we have for $\lambda$ sufficiently negative, $$
\aligned
0&\leq F_{\beta}(v_\lambda(\bar x)) -F_{\beta}(v(\bar x)) +c_2(\bar x)U(\bar x)\\
&\leq \frac{CV(\bar x)}{|\bar x|^\beta}+c_2(\bar x)U(\bar x)\\
&\leq C(\frac{V(\bar x)}{|\bar x|^\beta}-c_2(\bar x)c_1(\tilde x)|\tilde x|^\alpha V(\bar x))\\
&\leq\frac{CV(\bar x)}{|\bar x|^\beta}(1-c_1(\tilde x)|\tilde x|^\alpha c_2(\bar x)|\bar  x|^\beta).
\endaligned
$$
It follows that $1\leq c_1(\tilde x)|\tilde x|^\alpha c_2(\bar x)|\bar  x|^\beta.$
However,  from  \eqref{eq:aa11} we have $ c_1(\tilde x)|\tilde x|^\alpha c_2(\bar x)|\bar  x|^\beta<1 $ for  $|\tilde x|$ and $|\bar x|$ sufficiently large.
This contradiction shows that \eqref{eq:aaa11} must be true.

\section{Symmetry of solutions in the whole space $R^n$ }

{\bf Proof of Theorem \ref{thma3}.}
Choose an arbitrary direction as the $x_1$-axis. Let $T_\lambda =\{x\in R^n \mid x_1=\lambda\},~x^\lambda=(2\lambda-x_1,x'),~ u_\lambda(x)=u(x^\lambda),~\Sigma_\lambda=\{x\in R^n|x_1<\lambda\}$,
$$U_\lambda(x)=u_\lambda(x)-u(x),~V_\lambda(x)=v_\lambda(x)-v(x).$$
Step1.  \textit{Start moving the plane $T_\lambda$ from $-\infty$ to the right in $x_1$-direction.}

We will show that for $\lambda$ sufficiently negative,
\begin{equation}\label{eq:a16}U_\lambda(x)\geq0,~V_\lambda(x)\geq0,~x\in \Sigma_\lambda.\end{equation}
For the fixed $\lambda$  and $x\in\Sigma_\lambda,$   by \eqref{eq:aa1},
$$u(x)\rightarrow0,~\mbox{as}~|x|\rightarrow+\infty.$$
As $|x|\rightarrow+\infty$, we have $|x^\lambda|\rightarrow+\infty$;  it follows that
$$u_\lambda(x)=u(x^\lambda)\rightarrow0.$$
Thus for $x\in \Sigma_\lambda  ,$ \begin{equation}\label{eq:aa17}U_\lambda(x)\rightarrow 0,~\mbox{as}~|x|\rightarrow+\infty.\end{equation}
Similarly, one can show that for $x\in \Sigma_\lambda  ,$ $$V_\lambda(x)\rightarrow 0,~\mbox{as}~|x|\rightarrow+\infty.$$
By the mean value theorem it is easy to see that
$$F_\alpha (u_\lambda(x))-F_\alpha (u(x))=f(v_\lambda(x))-f(v(x))=f'(\xi_\lambda(x))V_\lambda(x),$$
and $$F_\beta (v_\lambda(x))-F_\beta (v(x))=g(u_\lambda(x))-g(u(x))=g'(\eta_\lambda (x))U_\lambda(x),$$
where $\xi_\lambda(x)$ is valued between $v_\lambda(x)$ and $v(x)$;  $\eta_\lambda(x)$ is valued between $u_\lambda(x)$ and $u(x)$. By Theorem \ref{thma1}, it suffices to check the decay rate of $f'(\xi_\lambda(x))$ and $g'(\eta_\lambda(x))$,
  at the points   where $V_\lambda(x)$ and $U_\lambda(x)$ are negative respectively.
Since $u_\lambda( x)<u(x)$ and $v_\lambda(x)<v(  x)$,
 we have $$0\leq u_\lambda( x)\leq \eta_\lambda(x)\leq u( x),~0\leq v_\lambda( x)\leq \xi_\lambda(x)\leq v( x).$$
At those points  for $|x|$ sufficiently large, the decay assumptions  \eqref{eq:aa1} and  \eqref{eq:aa2} instantly yields that
$$c_1(x)=f'(\xi_\lambda(x))\sim o(\frac{1}{| x|^\alpha}), ~ c_2(x)=g'(\eta_\lambda(x))\sim o(\frac{1}{| x|^\beta}).$$
Consequently, there exists $R_0>0,$  such that, if $\tilde x$ and $\bar x$ are negative minima of $U_\lambda(x)$ and  $V_\lambda(x)$ in $\Sigma_\lambda$ respectively, then by Theorem \ref{thma1},  it holds that
\begin{equation}\label{eq:aa16}|\tilde x|\leq R_0 ~\mbox{or}~|\bar x|\leq R_0.\end{equation}
Without loss of generality, we may assume
\begin{equation}\label{eq:aa18}|\tilde x|\leq R_0 .\end{equation} For $\lambda$ sufficiently negative, combining  \eqref{eq:aa17} with fact that $$U_\lambda(x)=0,~x\in T_\lambda,$$  we know if $U_\lambda(x)<0$ in $\Sigma_\lambda$, then $U_\lambda(x)$
must have a negative minimum in $\Sigma_\lambda$. This contradicts  \eqref{eq:aa18}. Hence, for $\lambda$ sufficiently negative we have \begin{equation}\label{eq:aa108}U_\lambda(x)\geq0,\end{equation} it follows that $V_\lambda(x)\geq 0$ in $\Sigma_\lambda$. Otherwise, there exists $\bar x$ in $\Sigma_\lambda$ such that $$V_\lambda(\bar x)=\underset{\Sigma_\lambda}{\min}V_\lambda(x)<0,$$
from \eqref{eq:aa13},  we have \begin{equation}\label{eq:aaa13}F_\beta(v_\lambda(\bar x))-F_\beta(v(\bar x))< 0.\end{equation}
However, combining   \eqref{eq:a11} with   \eqref{eq:aa108}, we have $F_\beta(v_\lambda(\bar x))-F_\beta(v(\bar x))\geq 0.$  This is a contradiction with \eqref{eq:aaa13} and $V_\lambda(x)$ cannot attain its negative value in  $\Sigma_\lambda.$
It follows that \eqref{eq:a16} must be true. This completes the preparation for the moving planes.

Step 2. \textit{Keep moving  the planes to the right to the limiting position $T_{\lambda_0}$ as long as
\eqref{eq:a16} holds}.

Let $$\lambda_0=\sup\{\lambda \mid  U_\mu(x),~ V_\mu(x) \geq 0, ~x\in \Sigma_\mu, ~\mu\leq \lambda\}.$$

Obviously, \begin{equation}\label{eq:a61}\lambda_0<\infty.\end{equation}
Otherwise, for any $\lambda>0,$
$$u(0^\lambda)>u(0)>0,~~v(0^\lambda)>v(0)>0.$$
Meanwhile, $$u(0^\lambda)\sim \frac{1}{|0^\lambda|^\beta},~~v(0^\lambda)\sim \frac{1}{|0^\lambda|^\alpha},~\lambda\rightarrow\infty.$$
This is a contradiction and \eqref{eq:a61} is proved.

Now,  we point out that \begin{equation}\label{eq:a18} U_{\lambda_0}(x)\equiv0, ~ V_{\lambda_0}(x) \equiv 0,~ x\in \Sigma_{\lambda_0}.\end{equation}
Otherwise,  we will show that the plane $T_\lambda$ can be moved further to the right. More rigorously,  there exists some $\epsilon>0,$ such that for any
$\lambda\in [\lambda_0,\lambda_0+\epsilon)$ we have
\begin{equation}\label{eq:a17}U_\lambda(x)\geq 0,~V_\lambda(x)\geq 0, ~x\in \Sigma_\lambda.\end{equation}
This is a contradiction with the definition of $\lambda_0$. Hence we must have \eqref{eq:a18}.

Now we to prove \eqref{eq:a17} by using Theorem \ref{thma2} and Theorem \ref{thma1}.

Suppose \eqref{eq:a18} is false, then $U_{\lambda_0}(x)\geq 0$ and $V_{\lambda_0}(x)\geq 0$ are positive somewhere in $\Sigma_{\lambda_0}$, and Theorem \ref{thma2} gives $$U_{\lambda_0}(x)> 0,~V_{\lambda_0}(x)> 0, ~x\in \Sigma_{\lambda_0} .$$
Let $R_0$ be determined in Theorem \ref{thma1}. It follows that for any  $\delta>\epsilon>0$,
$$U_{\lambda_0}(x)\geq c_0>0,~V_{\lambda_0}(x)\geq c_0>0,~  x\in \overline{\Sigma_{\lambda_0-\delta}\cap B_{R_0}(0)} .$$
From the continuity of $U_\lambda(x)$ and $ V_\lambda(x)$ with respect to $\lambda$, there exists   $\epsilon>0,$  such that  for all  $\lambda\in [\lambda_0,\lambda_0+\epsilon),$  we have

\begin{equation}\label{eq:a19}
U_\lambda(x)\geq 0,~V_\lambda(x)\geq 0, ~x\in \overline{\Sigma_{\lambda_0-\delta}\cap B_{R_0}(0)} .
\end{equation}
Suppose that \eqref{eq:a17} is false,  we have $U_\lambda(x)<0,~V_\lambda(x)<0,~x\in \Sigma_\lambda$.  If $\tilde x$ and $\bar x$ are negative minima of  $U_\lambda(x)$ and $ V_\lambda(x)$ in $\Sigma_\lambda$ respectively. Next we  consider two possibilities.

\textbf{Case 1}. One of the negative minima of $U_\lambda(x)$ and $V_\lambda(x)$ lies in $B_{R_0}(0)$, i.e. it is in the narrow region $\Sigma_{\lambda_0+\epsilon}\backslash \Sigma_{\lambda_0-\delta}$. The other is outside of $B_{R_0}(0)$.
Without loss of generality, we may assume  the negative minimum of $U_\lambda(x)$ lies in  $B_{R_0}(0).$
 from \eqref{eq:a8}, we have
 \begin{equation}\label{eq:aa8}
U_\lambda(\tilde x)\geq -c_1(\tilde x)l^\alpha V_\lambda(\tilde x).
\end{equation}
Furthermore,  we know
$$\aligned
0&\leq F_\beta (v_\lambda(\bar x))-F_\beta (v(\bar x))+c_2(\bar x)U_\lambda(\bar x)\\
&\leq\frac{CV_\lambda(\bar x)}{|\bar x|^\beta}+c_2(\bar x)U_\lambda(\tilde x)\\
&\leq C\{\frac{V_\lambda(\bar x)}{|\bar x|^\beta}-c_2(\bar x)c_1(\tilde x)l^\alpha V_\lambda(\tilde x)\}\\
&\leq C\{\frac{V_\lambda(\bar x)}{|\bar x|^\beta}-c_2(\bar x)c_1(\tilde x)l^\alpha V_\lambda(\bar x)\}\\
&\leq C \frac{V_\lambda(\bar x)}{|\bar x|^\beta}[1-c_1(\tilde x)l^{\alpha}c_2(\bar x)|\bar x|^\beta].
\endaligned$$
Hence \begin{equation}\label{eq:aaa20}1\leq c_1(\tilde x)l^{\alpha}c_2(\bar x)|\bar x|^\beta.\end{equation}
From  \eqref{eq:aa11},  we know that $c_2(\bar x)|\bar x|^\beta$ is small for $|\bar x|$ sufficiently large. Since $l=\epsilon+\delta$ is very narrow and $c_1(\tilde x)$ is bounded from below in
$\Sigma_{\lambda_0+\epsilon}\backslash \Sigma_{\lambda_0-\delta}$, $c_1(\tilde x)l^{\alpha}$ can be small. Consequently,  $c_1(\tilde x)l^{\alpha}c_2(\bar x)|\bar x|^\beta<1$. This is a contradiction with \eqref{eq:aaa20} and \eqref{eq:a17} is proved.

\textbf{Case 2}. The negative minima of $U_\lambda(x)$ and $V_\lambda(x)$ lie in $B_{R_0}(0)$, i.e. they are all in the narrow region $\Sigma_{\lambda_0+\epsilon}\backslash \Sigma_{\lambda_0-\delta}$.

By \eqref{eq:a7},
\begin{equation}\label{eq:a23} F_{\alpha}(u_\lambda(\tilde x)) -F_{\alpha}(u(\tilde x)) \leq \frac{CU_\lambda(\tilde x)}{l^\alpha}<0,
\end{equation}  where $l=\delta+\epsilon.$
Together  with \eqref{eq:a3}, it implies
\begin{equation}\label{eq:aa8}
U_\lambda(\tilde x)\geq -c_1(\tilde x)l^\alpha V_\lambda(\tilde x).
\end{equation}
Similarly to \eqref{eq:a23}, we  derive
$$F_\beta (v_\lambda(\bar x))-F_\beta (v(\bar x))\leq \frac{C V_\lambda(\bar  x)}{l^\beta}<0.$$
Noting \eqref{eq:aa8}, we have for $l$ sufficiently small,
$$\aligned
0&\leq F_\beta (v_\lambda(\bar x))-F_\beta (v(\bar x))+c_2(\bar x)U_\lambda(\bar x)\\
&\leq\frac{CV_\lambda(\bar x)}{l^\beta}+c_2(\bar x)U_\lambda(\tilde x)\\
&\leq C\{\frac{V_\lambda(\bar x)}{l^\beta}-c_2(\bar x)c_1(\tilde x)l^\alpha V_\lambda(\tilde x)\}\\
&\leq C\{\frac{V_\lambda(\bar x)}{l^\beta}-c_2(\bar x)c_1(\tilde x)l^\alpha V_\lambda(\bar x)\}\\
&\leq C \frac{V_\lambda(\bar x)}{l^\beta}[1-c_1(\tilde x)c_2(\bar x)l^{\alpha+\beta}]\\
&<0.
\endaligned$$
This contradiction shows that \eqref{eq:a17} must be true.

Now we have shown that $ U_{\lambda_0}(x)\equiv 0,~V_{\lambda_0}(x)\equiv0, ~ x\in \Sigma_{\lambda_0}.$
Since the $x_1$ direction can be chosen arbitrarily, we  actually prove that $u(x)$ and $v(x)$ must be radially symmetric  about some point $x^0.$
Also the monotonicity follows easily from the argument.

This completes the proof of  Theorem \ref{thma3}.
\section{Non-existence of solutions on a half space $R^n_+$ }
We investigate  the system \eqref{eq:a2}.

{\bf Proof of Theorem \ref{thma4}.} Based on \eqref{eq:a28} and $f(0)=0,~g(0)=0$, one can see
from the proof of Lemma \ref{lemaa1} that
$$\mbox{either} ~u(x)>0,~v(x)>0~ \mbox{or} ~u(x)\equiv 0,~v(x)\equiv 0 , ~\mbox{ for}~ x\in R^n_+.$$
In fact, without loss of generality, assume $u(x)\not \equiv 0,$ there exists  $x^0$ such that $u(x^0)=0$, and
$$F_\alpha(u(x^0))=c_{n,\alpha} PV \int _{R^n}\frac{G(u(x^0)-u(y))}{|x^0-y|^{n+\alpha}}dy<0,$$ i.e. $0\leq f(v(x))= F_\alpha(u(x))<0 ,$  this is impossible.
 Hence if $u(x)$ or  $v(x)$
 attains 0 somewhere in $R^n_+$, then $u(x)=v(x)\equiv 0, x\in R^n_+.$

Hence in the following, we  assume that $u(x)>0$ and $v(x)>0$ in $R^n_+.$
Let us carry on the method of moving planes on the solution $u$ along $x_n$ direction.

Dote $T_\lambda=\{x\in R^n| x_n=\lambda\},~\lambda>0,$ ~$\Sigma_\lambda=\{x\in R^n| 0<x_n<\lambda\}.$
Let $x^\lambda=(x_1,\cdots,x_{n-1},2\lambda-x_n)$ be the reflection of $x$ about the plane $T_\lambda$, and $U_\lambda(x)=u_\lambda(x)-u(x),~V_\lambda(x)=v_\lambda(x)-v(x).$

The key ingredient \eqref{eq:aa20} is obtained in this proof of Theorem \ref{thma2}. To see that it still applies in this situation, we only need to take $\Sigma=\Sigma_\lambda\cup R^n_-,$
where $R^n_-=\{x\in R^n| x_n\leq 0\}.$

Step1. For $\lambda$ sufficiently small, we have immediately
\begin{equation}\label{eq:a26}
U_\lambda(x)\geq 0, V_\lambda(x)\geq 0,~x\in \Sigma_\lambda,\end{equation}
since $\Sigma_\lambda$ is a narrow region.

Step2. Since \eqref{eq:a26} provides a starting point,  we  move the plane $T_\lambda$ upward as long as \eqref{eq:a26} holds.
Define $$\lambda_0=\sup\{\lambda>0|U_\mu(x)\geq 0, V_\mu(x)\geq 0,~x\in \Sigma_\mu,~ \mu \leq \lambda \}.$$
We show that \begin{equation}\label{eq:a27}\lambda_0=\infty.\end{equation}
Otherwise, if $\lambda_0<\infty,$ then using \eqref{eq:a26},  Theorem \ref{thma2},  Theorem \ref{thma1} and going through
the similar arguments as in  Section 3, we  are able to show
$$U_{\lambda_0}\equiv 0,~ V_{\lambda_0}\equiv 0, ~x\in \Sigma_{\lambda_0}, $$
which implies
$$u(x_1,\cdots, x_{n-1},2\lambda_0)=u(x_1,\cdots, x_{n-1},0)=0,$$
$$v(x_1,\cdots, x_{n-1},2\lambda_0)=v(x_1,\cdots, x_{n-1},0)=0.$$
This is impossible, because we assume that $u(x),v(x)>0$ in $R^n_+.$

Therefore, \eqref{eq:a27} must be valid  and the solutions $u(x),v(x)$
are  increasing with respect to $x_n.$ This contradicts \eqref{eq:a28}.

This completes the proof of  Theorem \ref{thma4}.

\end{document}